\documentclass[12pt,a4paper]{amsart}

\usepackage[english]{babel}
\usepackage{lmodern}
\usepackage{newlfont}
\usepackage{booktabs}
\usepackage{color}
\usepackage[T1]{fontenc}
\usepackage[utf8]{inputenc}
\usepackage{indentfirst}
\usepackage{amsmath}
\usepackage{amsmath,amssymb}
\usepackage{amsthm}
\usepackage{geometry}
\usepackage{mathtools}
\usepackage{listings}
\usepackage{amsfonts}
\usepackage{braket}
\usepackage{emptypage}
\usepackage{newlfont}
\usepackage{amssymb}
\usepackage{graphicx}
\usepackage[italian]{varioref}
\usepackage{accents}
\usepackage{comment}
\usepackage{stmaryrd}
\usepackage{pgfplots}
\makeatletter
\labelformat{equation}{\tagform@{#1}}
\makeatother
\usepackage{hyperref}
\usepackage{relsize}
\usepackage{faktor}
\usepackage{enumerate}
\usepackage{fancyhdr}
\usepackage[colorinlistoftodos]{todonotes}

\renewcommand{\phi}{\varphi}
\renewcommand{\epsilon}{\varepsilon}
\newcommand{\al}{\alpha}
\newcommand{\numberset}{\mathbb}

\newcommand{\R}{\numberset{R}}

\newcommand{\la}{\lambda}
\newcommand{\Om}{\Omega}

\theoremstyle{definition}

\theoremstyle{definition}                                                                         
\newtheorem{definizione}{Definizione}[section]
\theoremstyle{definition}                                                                         

\theoremstyle{plain}                                                                              
\newtheorem{thm}[definizione]{Theorem}
\theoremstyle{plain}    
\newtheorem{prop}[definizione]{Proposition}
\theoremstyle{plain}     
\newtheorem{lemma}[definizione]{Lemma}
\theoremstyle{plain}
\newtheorem{cor}[definizione]{Corollary}
\theoremstyle{definition}                                                                         
\newtheorem{es}{Example}
\newtheorem{remark}{Remark}

\numberwithin{equation}{section}
\begin{document}
\title[On the linearization.. ]{On the linearization of the supercritical Lane--Emden equation in $\mathbb{R}^N$}

\author[Grossi]{Massimo Grossi }
\address{Dipartimento di Scienze di Base e Applicate per l'Ingegneria, Universit\`a di Roma ``La Sapienza'', Via Scarpa 12 - 00161 Roma, Italy, e-mail: {\sf massimo.grossi@uniroma1.it}.}
\author[Neves]{Sérgio L. N Neves}
\address{Departamento de Matemática, Universidade Estadual Paulista, Rua Cristóvão Colombo 2265, São José do Rio Preto, SP, 15054-000, Brazil, e-mail: {\sf sergio.neves@unesp.br}.}
\author[Wynants]{Zoe Saudade Wynants}
\address{Dipartimento di Scienze di Base e Applicate per l'Ingegneria, Universit\`a di Roma ``La Sapienza'', Via Scarpa 12 - 00161 Roma, Italy, e-mail: {\sf zoesaudade.wynants@uniroma1.it}.}

\maketitle

\begin{abstract}
Let $u$ be a positive radial solution of the supercritical
Lane-Emden equation
\[
    -\Delta u=u^p \qquad \text{in }\mathbb{R}^N,
    \qquad
    p>\frac{N+2}{N-2}.
\]
In this paper we classify the bounded solutions of  the associated linearized operator $-\Delta w -pu^{p-1}w=0$ in $\R^N$.
The resulting classification reveals the role of the threshold value
$p=\frac{N + \sqrt{N^2 - 8N}}{4}$.
In particular, above this threshold, the solutions are generated by the invariance of the problem under dilations and translations. Below this threshold, many new solutions appear, and their number tends to infinity as $p$ approaches the Joseph–Lundgren exponent and $N$ is large.
\end{abstract}

\section{Introduction}

We consider positive entire solutions of the Lane--Emden equation
\begin{equation}\label{eq:LE}
    -\Delta u=u^p \qquad \text{in }\mathbb R^N,
\end{equation}
where $N\geq 3$ and
\[
    p> p_S:=\frac{N+2}{N-2}.
\]
In the supercritical regime equation \ref{eq:LE} admits a family of positive radial
entire solutions, which are uniquely determined by their
value at the origin. More precisely, for every $a>0$,
there exists a unique positive radial solution $u_a$
satisfying
\[
u_a(0)=a, \qquad u_a'(0)=0.
\]
These solutions are smooth, strictly decreasing in the
radial variable, and scaling invariant in the sense that if $u$ solves \ref{eq:LE} then
\begin{equation}\label{m1}
u_\sigma(r)=\sigma^{\frac{2}{p-1}} u\left(\sigma r\right)
\end{equation}
also satisfies \ref{eq:LE} for any $\sigma>0$.

In addition to this family of regular solutions, there exists a unique
positive radial solution which is singular at the origin. It is given
explicitly by
\begin{equation}\label{eq:singular-solution}
    u_s(r)=Lr^{-\frac{2}{p-1}},
\end{equation}
where
\begin{equation}\label{eq:L-constant}
    L^{p-1}
    =
    \frac{2}{p-1}
    \left(
        N-2-\frac{2}{p-1}
    \right).
\end{equation}
The singular solution describes the leading-order behavior of every
regular positive radial solution at infinity. Namely,
\begin{equation}\label{eq:regular-singular-asymptotic}
    \lim_{r\to+\infty}
    r^\frac2{p-1}u_a(r)=L\quad\hbox{for any }a>0.
\end{equation}
Thus, although the solutions $u_a$ are regular at the origin, they
approach at infinity the same homogeneous profile determined by the
singular solution $u_s$.
These properties are well known and follow from the
classical analysis of radial solutions of the
Lane--Emden equation; see, for instance, 
\cite{GuiNiWang1992, SJ}.
The purpose of this paper is to study the bounded kernel of the
linearized operator
\begin{equation}\label{eq:linearized_operator}
    \mathcal L_u w=-\Delta w-pu^{p-1}w
\end{equation}
where $u$ is a positive regular radial solution of \ref{eq:LE}.
For the sake of simplicity in the following, we will consider 
\begin{equation}\label{u}
u:=u_1
\end{equation}
hence
$$u(0)=1\quad\text{and} \lim_{r\to+\infty}
    r^\frac2{p-1}u(r)=L.$$
More precisely, we
are concerned with the classification of solutions
\begin{equation}\label{eq:linearized}
    \mathcal L_u w=0 \qquad \text{in }\mathbb R^N,
    \qquad w\in L^\infty(\mathbb R^N).
\end{equation}
Nondegeneracy properties of this type are a basic ingredient in the
analysis of nonlinear elliptic equations, in perturbative and
gluing constructions, blow-up analysis, and Lyapunov--Schmidt
reductions.

At the critical Sobolev exponent $p=p_S$, the corresponding question
is by now classical. In this case all positive solutions of
\ref{eq:LE} are, up to translations and dilations, given by the
Aubin--Talenti family
\[
    U_{\lambda,\xi}(x)
    =
    c_N
    \left(
        \frac{\lambda}
        {1+\lambda^2|x-\xi|^2}
    \right)^{\frac{N-2}{2}};
\]
see \cite{Aubin1976,Talenti1976,CaffarelliGidasSpruck1989}. The
linearized operator at a bubble $U_{\lambda,\xi}$ is nondegenerate modulo the natural
symmetries of the equation: its bounded kernel is generated by the
derivatives of the family $U_{\lambda,\xi}$ with respect to the
scaling and translation parameters. In other words, see \cite{Rey1990}, 
\begin{equation}\label{eq:critical_kernel}
    \ker_{L^\infty(\mathbb R^N)} \mathcal L_U
    =
    \operatorname{span}
    \left\{
        \frac{N-2}{2}U+x\cdot\nabla U,\,
        \partial_{x_1}U,\ldots,\partial_{x_N}U
    \right\}.
\end{equation}
Related
nondegeneracy results for critical Lane--Emden systems have also been
obtained in \cite{FrankKimPistoia2021}.

The supercritical case $p>p_S$ presents a substantially different
picture. Positive radial entire solutions still exist, but their
asymptotic behavior is governed by the singular radial solution, see \cite{GuiNiWang1992, SJ}. Consequently,
\begin{equation}\label{eq:inverse_square}
    pu^{p-1}(r)
    =
    \frac{pL^{p-1}}{r^2}(1+o(1))\hspace{3mm}\text{when}\hspace{3mm}r\to+\infty.
\end{equation}
Thus, unlike the critical case, where the linearized potential decays
as $|x|^{-4}$, in the supercritical regime the potential has the
critical inverse-square decay. This observation accounts for a substantial change in the spectral properties of the
linearized operator.

The inverse-square structure is also at the heart of the stability
theory for supercritical Lane--Emden equations. Starting with the
classical work of Joseph and Lundgren \cite{JosephLundgren1973}, it
has been understood that the behavior of radial solutions is governed
by a second critical exponent, now usually referred to as the
Joseph--Lundgren exponent
\[
p_{\mathrm{JL}}(N):=
\begin{cases}
+\infty, & 2\leq N\leq 10,\\[6pt]
\displaystyle
\frac{(N-2)^2-4N+8\sqrt{N-1}}
     {(N-2)(N-10)},
& N\geq 11.
\end{cases}
\]
It is now well established that, for \(p<p_\mathrm{JL}\), all nontrivial solutions to the Lane–Emden equation in \(\mathbb{R}^N\) are unstable (\cite{Farina2007}). This instability result applies, in particular, to the radial solution $u_a$.
Stability and finite Morse index properties
of entire solutions have subsequently been investigated in several
works; we refer in particular to
\cite{Farina2007,DupaigneFarina2010} and to the monograph
\cite{Dupaigne2011} for a comprehensive account.\\
The linear theory associated with supercritical Lane--Emden profiles
has previously been investigated in connection with perturbative
constructions for supercritical elliptic equations. In particular,
Dávila, del Pino, Musso and Wei \cite{DavilaDelPinoMussoWei2007}
developed a linear theory around positive radial solutions of
\[
    -\Delta u=u^p \qquad \text{in }\mathbb R^N,
\]
in their study of standing waves for supercritical nonlinear
Schrödinger equations. Their analysis is carried out in weighted
spaces adapted to the asymptotic behavior of the underlying
supercritical profile and to the perturbative scheme. A comparison with our result is provided in Remark \ref{confronto}.
Collot in \cite{Collot2017} also studies this linearized operator, in connection with type II blow-up for the energy-supercritical semilinear heat equation; there, $p$ is assumed to be an odd integer, $p = 2q+1$, and $p > p_\mathrm{JL}$.\\
The question considered here is of a different nature. Rather than
prescribing a weighted decay at infinity, we seek to characterize the
kernel of the linearized operator
under the sole assumption that the solution is bounded. More
precisely, we classify all the solutions satisfying \ref{eq:linearized}.\\
This distinction is particularly relevant in the supercritical
regime: our results require a delicate analysis of the behavior of the solutions at infinity.

Indeed, the invariances of
\ref{eq:LE} for dilations and translations immediately yield the bounded solutions
\begin{equation}\label{eq:natural_kernel}
    \frac{2}{p-1}u+ru',
    \qquad
    \frac{\partial u}{\partial x_i},
    \quad i=1,\ldots,N,
\end{equation}
of the linearized equation. It is therefore natural to ask whether
these functions exhaust the bounded kernel of $\mathcal L_u$. While in the
critical case the answer is precisely the classical nondegeneracy
statement \ref{eq:critical_kernel}, the slower decay
\ref{eq:inverse_square} makes the corresponding supercritical
question more delicate. In particular, the potential and the angular
part of the Laplacian have the same homogeneity at infinity, so that
higher spherical modes cannot be discarded solely on the basis of
their asymptotic order.\\
Our first main result is the following.
\begin{thm}\label{main}
When $N\geq12$ and $p > p_2=\frac{N + \sqrt{N^2 - 8N}}{4}$ the kernel of the linear operator $\mathcal{L}_u=-\Delta-pu^{p-1}$ in $C^2(\mathbb R^N)\cap L^\infty(\mathbb R^N)$ is
    $$\ker_{L^\infty(\mathbb R^N)} \mathcal L_u=\operatorname{span}\left\{\frac{2}{p-1}u+ru', \frac{\partial u}{\partial x_1},\dots,\frac{\partial u}{\partial x_N}\right\}.$$
    When $N=11$ the same results hold for $p>p_\mathrm{JL}$.
\end{thm}
\begin{remark}
    Setting $v_0=\frac{2}{p-1}u+ru'$, we prove (see Corollary \ref{pos}) 
    $$v_0>0\hspace{2mm}\text{for all}\hspace{2mm} r\geq 0\hspace{2mm}\text{ when}\quad p\geq p_\mathrm{JL}.$$
    The positivity of $v_0$ highlights a significant difference from the critical case $p=\frac{N+2}{N-2}$, where the corresponding eigenfunction changes sign. This is due to the stability of the solution $u$ in the supercritical case, a property that fails for the standard bubble in the critical case. Although the positivity of $v_0$ appears to be known among experts in the field, we were unable to find a proof in the literature. For this reason, we provide one in Corollary \ref{pos} in Section \ref{s3}, where we also relate the positivity of $v_0$ to the stability of solutions of the Gelfand problem in bounded domains.
\end{remark}
The condition $p > p_2$ is sharp in Theorem \ref{main}; indeed, for $p_\mathrm{JL}(N) < p \le p_2$, new bounded elements appear in the kernel (Theorem \ref{main2} and Remark \ref{r1} below). In order to state our result we need to introduce some notations. Recall
\begin{equation}\label{eq:L}
L^{p-1}=\frac{2}{p-1}\left(N-2-\frac{2}{p-1}\right),
\end{equation}
and for $j\ge 0$ we denote by
\[
\mu_j=j(j+N-2)
\]
the eigenvalues of the Laplace--Beltrami operator $-\Delta_{\mathbb S^{N-1}}$, and by
$\mathcal H_j$ the associated eigenspace, namely the space of spherical harmonics
of degree $j$. We recall that
\begin{equation}\label{eq:dimHj}
d_j:=\dim\mathcal H_j=\binom{N+j-1}{j}-\binom{N+j-3}{j-2},
\end{equation}
with the convention $\binom{n}{k}=0$ for $k<0$.
Finally, we set
\begin{equation}\label{eq:J-Np-JL}
J(N,p):=\frac{\sqrt{(N-2)^2+4pL^{p-1}}-(N-2)}{2},
\end{equation}
so that, for every integer $j\ge 0$, $\mu_j\le pL^{p-1}$ if and only if $j\le J(N,p)$.

\begin{thm}\label{main2}
Let $N\ge 12$ and $p_{\mathrm{JL}}(N)<p\le p_2=\frac{N + \sqrt{N^2 - 8N}}4$.
Then the kernel of  $\mathcal L_u$ is given by
\begin{equation}\label{kerfin}
\begin{split}
    \ker_{L^\infty(\mathbb R^N)} \mathcal L_u=\operatorname{span}&\left\{\frac{2}{p-1}u+ru', \frac{\partial u}{\partial x_1},\dots,\frac{\partial u}{\partial x_N}\right\}\bigoplus_{j=2}^{\lfloor J(N,p)\rfloor}\left\{v_j(|x|)\,\varphi\big(\tfrac{x}{|x|}\big):\ \varphi\in\mathcal H_j\right\}.
    \end{split}
\end{equation}
Here $v_j$ is the positive solution of 
\begin{equation}\label{eq:radial}
v''+\frac{N-1}{r}v'+\left(pu(r)^{p-1}-\frac{\mu_j}{r^2}\right)v=0,
\quad\text{for}\quad r>0
\end{equation}
 which is regular at the origin and satisfies
\begin{equation}\label{rp}
v_j(r)\sim r^{\frac{-(N-2)+\sqrt{(N-2)^2+4(\mu_j-pL^{p-1})}}{2}}\text{ for }r\to+\infty.
\end{equation}
 In particular,
\begin{equation}\label{eq:kernel-dim-JL}
\dim\ker_{L^\infty(\mathbb R^N)}\mathcal L_u
=N+1+\sum_{j=2}^{\lfloor J(N,p)\rfloor}
\left[\binom{N+j-1}{j}-\binom{N+j-3}{j-2}\right],
\end{equation}
where each term of the sum is the dimension of the space of spherical harmonics of degree $j$.

Finally for every integer $m\ge 2$ there exist an integer $N\ge 12$ and an exponent $p>p_{\mathrm{JL}}(N)$ such that $\lfloor J(N,p)\rfloor=m$, that is, the sum in \ref{eq:kernel-dim-JL} consists of exactly $m-1$ terms.
\end{thm}
\begin{remark}\label{r1}
The sum in \ref{eq:kernel-dim-JL}  is nonempty if and only if $p\le p_2$; for $p>p_2$ it is empty and \ref{eq:kernel-dim-JL} reduces to Theorem~\ref{main}.
\end{remark}
\begin{remark}
 We have that $ J(N,p)\to \frac{(\sqrt2-1)(N-2)}{2}$ as $p\to p_\mathrm{JL}(N)$. So the number of such ``new solutions'' can be made arbitrarily large by taking $N$ large.
 \end{remark}
\begin{es}
For \(N=22\), the dimension of the bounded kernel increases significantly as \(p\) decreases towards the Joseph–Lundgren exponent. More precisely, one has
\[
\dim\ker_{L^\infty(\mathbb R^{22})}\mathcal{L}_u=
\begin{cases}
23, & p>\hat p_2,\hbox{ (standard dilations and translations)}\\[2mm]
275, & \hat p_3<p\le \hat p_2,\\[2mm]
2277, & \hat p_4<p\le \hat p_3,\\[2mm]
14674, & p_{\mathrm{JL}}(22)<p\le \hat p_4,
\end{cases}
\]
where
\[
p_{\mathrm{JL}}(22)
=1+\frac{4}{18-2\sqrt{21}}
\approx 1.4528,
\]
and
\[
\hat p_4=\frac{3}{2},\quad \hat p_3= 1+\frac{4}{18-4\sqrt{13}}
\approx 2.118,
\quad
\hat p_2=1+\frac{4}{18-2\sqrt{77}}
\approx 9.887.
\]
In particular, this example shows that the bounded kernel can have very large dimension even in the supercritical regime \(p>p_{\mathrm{JL}}(22)\).\\
Here $\hat p_j = \hat p_j(N)$ denotes the exponent defined in \ref{eq:pj-upper-JL} below,
namely the value of $p$ at which $\mu_j = pL^{p-1}$. In particular, $\hat p_2(N)$
coincides with the threshold $p_2$ of Theorem~\ref{main}.
\end{es}

The paper is organized as follows. In Section~2 we collect some preliminary results on the radial linearized equation, establishing existence, uniqueness and asymptotic behavior of its solutions. Section~3 is devoted to the proof of Theorem~\ref{main}, which shows that for $p>p_2$ the bounded kernel of $\mathcal L_u$ is spanned by the scaling and translation derivatives. In Section~4 we prove Theorem~\ref{main2}, characterizing the additional solutions that appear for $p\le p_2$ and determining the exact dimension of the kernel. Finally, in Section~5 we prove the positivity of $v_0$ via a blow-up argument, which relates the Lane--Emden equation to a Gelfand-type problem.

\section{Preliminaries}
Let us start this section with a lemma that will be useful in what follows.
\begin{lemma}\label{lem:existence-vj}
Given $u$ as in \ref{u}, for any $\mu\geq0$ the equation
\begin{equation}\label{eq:radial}
v''+\frac{N-1}{r}v'+\left(pu(r)^{p-1}-\frac{\mu}{r^2}\right)v=0,
\quad\text{for}\quad r>0,
\end{equation}
admits a unique nontrivial solution for all $r>0$ satisfying
\begin{equation}\label{eq:normalization}
\lim_{r\downarrow0}\frac{v(r)}{r^\alpha}=1
\end{equation}
where $\alpha$ is the nonnegative root of $\mu=\alpha(\alpha+N-2)$, namely
$$\alpha
=\frac{-N+2+\sqrt{(N-2)^2
+4\mu
}}{2}$$
\end{lemma}
\begin{proof}
The substitution
$v(r)=r^{\alpha}w(r)$ reduces \ref{eq:radial} to
\[
w''+\frac{N-1+2\alpha}{r}w'+pu(r)^{p-1}w=0.
\]
Put $d=N-1+2\alpha> -1$. The regular initial conditions
$w(0)=1$ and $w'(0)=0$ are equivalent to the Volterra equation
\begin{equation}\label{eq:volterra}
w(r)=1-p\int_0^r t^{-d}
\left(\int_0^t s^{d}u(s)^{p-1}w(s)\,ds\right)dt.
\end{equation}
Because $u$ is smooth near $0$ and $d>-1$, the integral operator in
\ref{eq:volterra} is a contraction on $C([0,R])$ for sufficiently small $R>0$. It has a unique solution $w$ on $[0,R]$.
In fact, let us fix a small \(R>0\) and $X=C([0,R])$. In this space define 
\[
(Tw)(r)
=
1-p\int_0^r t^{-d}
\left(
\int_0^t
s^{d}u(s)^{p-1}w(s)\,ds
\right)dt.
\]
A solution of \ref{eq:volterra} is a fixed point of $T$.
Since \(u\) is continuous near the origin, there exists  \(M>0\) such that
\[
u(s)^{p-1}\le M
\qquad (0\le s\le R).
\]
If \(w_1,w_2\in X\) we have (using that $d>-1$),
\[
|Tw_1(r)-Tw_2(r)|
\le{}
pM\|w_1-w_2\|_\infty\int_0^r t^{-d}
\left(\int_0^t s^{d}\,ds\right)dt\le 
\frac{pMR^2}{2(d+1)}\|w_1-w_2\|_\infty
\]
and then taking \(R\) sufficiently small we have that \(T\) is a contraction. By the Banach fixed-point theorem \(T\) has a unique fixed point \(w\in C([0,R])\).

Moreover it is not difficult to verify that the integral equation gives
$w(r)=1+O(r^2)$ and $w'(r)=O(r)$ as $r\downarrow0$.\\

Thus $v=r^{\alpha}w$ solves \ref{eq:radial} and satisfies
\ref{eq:normalization}, so $v\not\equiv0$.
Finally, the coefficients of \ref{eq:radial} are continuous on every compact
subinterval of $(0,\infty)$. Standard continuation for linear ordinary
differential equations extends $v$ uniquely to all $r>0$.
\end{proof}

We now describe the behavior of solutions of \ref{eq:radial} near the origin and at infinity.
\begin{lemma}\label{lem:asymptotics}
Let $N\ge 11$, $p>p_\mathrm{JL}(N)$ and $\mu\ge 0$, and set
\begin{equation}\label{eq:roots0}
\alpha^\pm=\frac{-(N-2)\pm\sqrt{(N-2)^2+4\mu}}{2},
\end{equation}
\begin{equation}\label{eq:roots1}
\beta^\pm=\frac{-(N-2)\pm\sqrt{(N-2)^2+4(\mu-pL^{p-1})}}{2}.
\end{equation}
Then $\alpha^\pm,\beta^\pm\in\R$ with $\alpha^-<0\le\alpha^+$ and $\beta^-<\beta^+$. Equation \ref{eq:radial} admits two linearly independent solutions $v^\pm$ such that
\begin{itemize}
\item[(i)]
Their asymptotic behaviors at the origin are
$$
v^\pm(r)\sim r^{\alpha^\pm}\quad\text{as } r\to 0^+ .
$$
In particular, up to a multiplicative constant, $v^+$ is the only solution bounded near the origin.
\item[(ii)] As $r\to+\infty$, either
$$v^+(r)\sim r^{\beta^+},\quad r\,(v^+)'(r)\sim\beta^+r^{\beta^+}\quad\text{or}\quad v^+(r)\sim r^{\beta^-},\quad r\,(v^+)'(r)\sim\beta^-r^{\beta^-}.$$
\item[(iii)] If $\mu\le pL^{p-1}$, $v^+$ is bounded at infinity. If $\mu>pL^{p-1}$, $v^+$ is bounded at infinity if and only if $v^+(r)\sim r^{\beta^-}$.
\end{itemize}
\end{lemma}

\begin{proof}
Multiplying \ref{eq:radial} by $r^2$, we obtain
\begin{equation}\label{eq:euler}
r^2v''(r)+(N-1)rv'(r)+\big(p\,u^{p-1}(r)\,r^2-\mu\big)v(r)=0 .
\end{equation}

\emph{Behavior near the origin.}
Since $u(0)=1$ and $u$ is a smooth radial solution of an equation with nonlinearity analytic near $1$, the function $p\,u^{p-1}(r)\,r^2$ is analytic in $r$ near $0$ and vanishes at $r=0$. Hence $r=0$ is a regular singular point of \ref{eq:euler}, and we may apply the Frobenius method (see \cite{BoyceDiPrima2005}). Looking for a solution of the form
$$
v(r)=r^{\alpha}\sum_{n=0}^{\infty}a_n r^n,\qquad a_0\neq 0,
$$
where the series converges in a neighborhood of the origin and substituting into \ref{eq:euler} we obtain
\begin{equation}\label{eq:series0}
\sum_{n=0}^{\infty}
\left[
(n+\alpha)(n+\alpha+N-2)-\mu
\right]a_nr^{n+\alpha}+p\sum_{n=0}^{\infty}u^{p-1}(r)a_nr^{n+\al+2}=0.
\end{equation}
When $r$ approaches $0$, since $u(0)=1$,
$$pu^{p-1}\to p$$
so we can rewrite \ref{eq:series0} as
\begin{equation}\label{eq:series}
\sum_{n=0}^{\infty}
\left[
(n+\alpha)(n+\alpha+N-2)-\mu
\right]a_nr^{n+\alpha}+{ \big(p+o(1)\big)}\sum_{n=0}^{\infty}a_nr^{n+\al+2}=0.
\end{equation}
The coefficient of the lowest power $r^{\alpha}$ must vanish. This gives the indicial equation
\begin{equation}\label{eq:indicial0}
\alpha(\alpha+N-2)-\mu=0,
\end{equation}
whose roots are $\alpha^\pm$. The Frobenius theory then yields two linearly independent solutions $v^\pm\sim r^{\alpha^\pm}$ as $r\to 0^+$; if $\alpha^+-\alpha^-$ is an integer, $v^-$ may contain a logarithmic term, which does not affect its leading behavior. Since $\alpha^-<0\le\alpha^+$, only $v^+$ is bounded near the origin, which proves (i).

\emph{Behavior at infinity.}
By \ref{eq:inverse_square}, $p\,u^{p-1}(r)\,r^2\to pL^{p-1}$ as $r\to+\infty$, so \ref{eq:euler} is asymptotically the Euler--Cauchy equation
$$
r^2v''+(N-1)rv'+\big(pL^{p-1}-\mu\big)v=0 .
$$
Looking for solutions of this limit equation of the form $v(r)=r^\beta$, we obtain the indicial equation
\begin{equation}\label{eq:indicialinf}
\beta(\beta+N-2)+pL^{p-1}-\mu=0,
\end{equation}
whose roots are $\beta^\pm$. Since $p>p_\mathrm{JL}$, we have $pL^{p-1}<\frac{(N-2)^2}{4}$, hence
$$
\frac{(N-2)^2}{4}+\mu-pL^{p-1}>0\qquad\text{for all }\mu\ge 0,
$$
and the roots $\beta^\pm$ are real and distinct.

To make this rigorous, we recall that $p\,u^{p-1}r^2\to pL^{p-1}$ at a polynomial rate (see \cite{GuiNiWang1992}). By standard asymptotic theory for linear ODEs, \ref{eq:euler} admits two linearly independent solutions behaving like $r^{\beta^+}$ and $r^{\beta^-}$ at infinity; since $v^+$ is a linear combination of them, (ii) follows.\\
Note that $\beta^-<-\frac{N-2}{2}<0$, and that $\beta^+>0$ if and only if $\mu>pL^{p-1}$. If $\mu\le pL^{p-1}$, both $r^{\beta_\pm}$ are bounded at infinity, while if $\mu>pL^{p-1}$ only $r^{\beta_-}$ is. This proves (iii).
\end{proof}

\section{The nondegeneracy case $p > p_2=\frac{N + \sqrt{N^2 - 8N}}{4}$}

We now prove Theorem \ref{main}
\begin{proof}
Let $u$ be given by \ref{u}. We are looking for bounded solutions of $\mathcal L_u w=0$.
Let $\mathcal H_j$ be the space of spherical harmonics of degree $j$, so that 
$$-\Delta_{\mathbb S^{N-1}}\varphi=\mu_j\varphi$$
for every $\varphi\in\mathcal H_j$, with $\mu_j=j(j+N-2)$. Expanding $w$ in spherical harmonics, each radial coefficient of degree $j$ solves
\begin{equation}\label{v}
v''+\frac{N-1}{r}v'+pu^{p-1}v-\frac{\mu_j}{r^2}v=0 ,
\end{equation}
and is bounded near the origin, since $w$ is bounded. By Lemma~\ref{lem:asymptotics}(i), it is therefore a multiple of the solution of \ref{eq:radial} given by Lemma \ref{lem:existence-vj}, which we denote by $v_j$. Hence
$$
w(r,\omega)=\sum_{j\ge 0} v_j(r)\,\varphi_j(\omega),\qquad \varphi_j\in\mathcal H_j ,
$$
and it suffices to determine for which $j$ the function $v_j$ is bounded at infinity.\\
We already have two solutions for \ref{v} which coincide, up to multiplicative constants, with the ones given by Lemma~\ref{lem:existence-vj}:
\begin{enumerate}
    \item when $j=0$, so $\mu_0=0$, then $v_0=\frac{2}{p-1}u+ru'$;
    \item when $j=1$, so $\mu_1=N-1$, then $v_1=u'$.
\end{enumerate}
When $p\geq p_\mathrm{JL}$ we already know from \cite{GuiNiWang1992} that
\begin{equation}\label{gui}
 v_0(0) = \frac{2}{p-1} \quad \text{and}\quad v_0\sim r^{\frac{-(N-2) + \sqrt{(N-2)^2 - 4pL^{p-1}}}{2}}\text{ for }r\to+\infty.
\end{equation}
Moreover
$$v_1(0)=0, \quad v_1(r)<0\quad\forall r>0\quad\text{and}\quad |v_1|\sim r^{-\frac{2}{p-1}-1}\text{ for }r\to+\infty.$$
By Lemma \ref{lem:asymptotics} 
$$v_j(r)\sim r^{\frac{-(N-2)+\sqrt{(N-2)^2+4\mu_j}}{2}}\text{ for }r\to0$$
and
$$v_j(r)\sim r^{\frac{-(N-2)\pm\sqrt{(N-2)^2+4(\mu_j-pL^{p-1})}}{2}}\text{ for }r\to+\infty.$$
Let us prove that $v_j$ is bounded if and only if $j\in\{0,1\}$. Assume, by contradiction, that $v_j$ is bounded for some $j\ge 2$.\\
Since $N\ge 12$, we have $p_2>p_{\mathrm{JL}}(N)$, and hence $v_0>0$ by Corollary~\ref{pos}. Moreover, the function $pL^{p-1}=\frac{2p}{p-1}\big(N-\frac{2p}{p-1}\big)$ is decreasing in $p$ for $p>p_{\mathrm{JL}}$ and equals $\mu_2$ at $p=p_2$. Since $\mu_j\ge\mu_2$ for $j\ge 2$, it follows that
\begin{equation}\label{pp2}
\mu_j>pL^{p-1}\qquad\text{for all } p>p_2 \text{ and } j\ge 2.
\end{equation}
As a consequence, when $p > p_2 = \frac{N + \sqrt{N^2 - 8N}}{4}$ the bounded solution $v_j$ of \ref{v} behaves like
$$v_j(r)\sim r^{\frac{-(N-2)-\sqrt{(N-2)^2+4(\mu_j-pL^{p-1})}}{2}}\text{ for }r\to+\infty$$
and
\begin{equation}\label{a0}
r^{N-3}v_0v_j\sim r^{-1-A}\text{ for }r\to+\infty
\end{equation}
with $A=\frac{\sqrt{(N-2)^2+4(\mu_j-pL^{p-1})}-\sqrt{(N-2)^2-4pL^{p-1}}}{2}>0$.\\
Using the equation satisfied by $v_0$ and $v_j$ we deduce that

\begin{equation}\label{a3}
-R^{N-1} \left( v_j' v_0 \right)(R) + R^{N-1} \left( v_0' v_j \right)(R) + \mu_j \int_0^R r^{N-3} v_jv_0 = 0
\end{equation}
By \ref{a0} and Lemma \ref{lem:asymptotics} (ii), the boundary terms in \ref{a3} tend to zero as $R \to +\infty$, hence
$$
\mu_j \int_0^{+\infty} r^{N-3} v_j v_0 = 0
$$
then, by  Corollary \ref{pos}, $v_j$ must change sign. Next we show that $v_j$ has a unique zero.
\\ By the Sturm Comparison Theorem, we can conclude that between any two zeros of $v_j$ there must be at least one zero of $v_0$. Therefore, if $v_j$ had more than one zero, then $v_0$ would have to vanish. Since $v_0>0$ everywhere by Corollary \ref{pos}, it follows that $v_j$ has exactly one zero.\\
Let $a$ denote the unique zero of $v_j$ and suppose that $v_j(r)>0$ if $r<a$, so $v_j'(a)\leq0$.\\
Taking $R=a$ in \ref{a3}
$$\underbrace{-a^{N-1} v_j'(a) v_0(a)}_{\geq0} + a^{N-1} v_0'(a) \underbrace{v_j(a)}_{=0} + \underbrace{\mu_j \int_0^a r^{N-3} v_jv_0}_{>0} = 0$$
thus we arrive at a contradiction.\\
When $N=11$ we have $p_2<p_\mathrm{JL}$, hence the same argument holds for all $p>p_\mathrm{JL}.$
\end{proof}
\begin{remark} \label{rem}
We point out that the hypothesis $p > p_2$ was used only to ensure that $\mu_j > pL^{p-1}$ for every $j \ge 2$. The argument above therefore shows, more generally, the following: if $N \ge 11$, $p > p_\mathrm{JL}(N)$ and $j \ge 2$ is such that $\mu_j > pL^{p-1}$, then \ref{v} has no nontrivial bounded solution.
\end{remark}

\section{Additional solutions for $p\leq p_2 =\frac{N + \sqrt{N^2 - 8N}}{4}$}
The next proposition contains some elementary computations concerning the function $ pL^{p-1}$.
\begin{prop}\label{prop:thresholds-above-JL}
Let $N\ge 11$ and $j\ge 2$ be integers. There exists an exponent $p>p_{\mathrm{JL}}(N)$ such that $\mu_j\le pL^{p-1}$ if and only if
\begin{equation}\label{eq:JL-dimensional-condition}
N>2+2(1+\sqrt2)\,j .
\end{equation}
In this case, setting
\begin{equation}\label{eq:pj-upper-JL}
\hat p_j(N):=1+\frac{4}{N-4-\sqrt{N^2-4j(j+N-2)}},
\end{equation}
we have $\hat p_j(N)>p_{\mathrm{JL}}(N)$ and, for $p>p_{\mathrm{JL}}(N)$,
\begin{equation}\label{eq:mode-interval-JL}
\mu_j\le pL^{p-1}\quad\Longleftrightarrow\quad p\le \hat p_j(N).
\end{equation}
If \ref{eq:JL-dimensional-condition} fails, then $\mu_j>pL^{p-1}$ for every $p>p_{\mathrm{JL}}(N)$.
\end{prop}

\begin{proof}
We have
\begin{equation}\label{eq:c-alpha-JL}
pL^{p-1}=\frac{N^2}{4}-\Big(\frac{2}{p-1}-\frac{N-4}{2}\Big)^2 .
\end{equation}
At the Joseph--Lundgren exponent,
$$
\frac{2}{p_{\mathrm{JL}}-1}=\frac{N-4-2\sqrt{N-1}}{2},\qquad
p_{\mathrm{JL}}L^{p_{\mathrm{JL}}-1}=\frac{N^2}{4}-(N-1)=\frac{(N-2)^2}{4}.
$$

\emph{Step 1.}
The condition $p>p_{\mathrm{JL}}$ is equivalent to 
$$0<\frac{2}{p-1}<\frac{2}{p_\mathrm{JL}-1}<\frac{N-4}{2}$$
Hence the function $p\mapsto pL^{p-1}$ is strictly decreasing on $(p_{\mathrm{JL}},+\infty)$ onto the open interval
$$
\Big(2(N-2),\,\frac{(N-2)^2}{4}\Big).
$$

\emph{Step 2.}
For $j\ge 2$ we have $\mu_j\ge\mu_2=2N>2(N-2)$. By Step 1, there exists $p>p_{\mathrm{JL}}$ with $\mu_j\le pL^{p-1}$ if and only if
\begin{equation}\label{eq:mu-cJL-strict}
\mu_j<\frac{(N-2)^2}{4},
\end{equation}
where the inequality is strict because the supremum $\frac{(N-2)^2}{4}$ is not attained. Substituting $\mu_j=j(j+N-2)$, \ref{eq:mu-cJL-strict} becomes
$$
(N-2)^2-4j(N-2)-4j^2>0,
$$
which is equivalent to \ref{eq:JL-dimensional-condition}. If \ref{eq:JL-dimensional-condition} fails, then $\mu_j\ge\frac{(N-2)^2}{4}>pL^{p-1}$ for every $p>p_{\mathrm{JL}}$, which proves the last assertion.

\emph{Step 3.}
Assume \ref{eq:JL-dimensional-condition}. By Step 1, $\mu_j=pL^{p-1}$ has exactly one solution $p>p_{\mathrm{JL}}$. By \ref{eq:c-alpha-JL}, this equation reads
$$
\left(\frac{2}{p-1}\right)^2-(N-4)\,\frac{2}{p-1}+\mu_j-2(N-2)=0,
$$
whose roots are
$$
\frac2{p^\pm-1}=\frac{N-4\pm\sqrt{N^2-4\mu_j}}{2}.
$$
The root 
$$\frac2{p^+-1}=\frac{N-4+\sqrt{N^2-4\mu_j}}{2}>\frac{N-4}{2}>\frac{2}{p_{\mathrm{JL}}-1}$$
so it does not correspond to an exponent $p>p_{\mathrm{JL}}$.\\
A direct computation shows that
$$
\frac2{p^--1}=\frac{N-4-\sqrt{N^2-4\mu_j}}{2}>0\iff\mu_j>2(N-2)$$
and
$$\frac{N-4-\sqrt{N^2-4\mu_j}}{2}<\frac{2}{p_{\mathrm{JL}}-1}\iff\mu_j<\frac{(N-2)^2}{4},
$$
and both conditions hold. The corresponding exponent is
$$
p^-=1+\frac{4}{N-4-\sqrt{N^2-4\mu_j}}=\hat p_j(N)>p_{\mathrm{JL}}(N).
$$
Finally, since $pL^{p-1}$ is strictly decreasing on $(p_{\mathrm{JL}},\infty)$ and equals $\mu_j$ at $p=\hat p_j(N)$, we obtain \ref{eq:mode-interval-JL}.
\end{proof}

We are now in position to prove Theorem \ref{main2}.
\begin{proof}[Proof of Theorem \ref{main2}]
The beginning of the proof is the same as in  Theorem \ref{main}. Using the decomposition in spherical harmonics we look for solutions $w$ of $\mathcal L_u(w)=0$ as
$$
w(r,\omega)=\sum_{j\ge 0} v_j(r)\,\varphi_j(\omega),\qquad \varphi_j\in\mathcal H_j.
$$
Of course, we again have the solutions corresponding to the indexes $j=0,1$, so let us consider the case $j\ge2$. 
Recall that for $p > p_\mathrm{JL}(N)$ and $j \ge 2$, the regular solution $v_j$ of \ref{v} is bounded at infinity if $\mu_j \le pL^{p-1}$, since in this case both indicial roots at infinity
$$\beta_j^\pm
=\frac{-N+2
\pm
\sqrt{
(N-2)^2
+4(\mu_j-pL^{p-1})
}}{2}$$
are nonpositive. Conversely, if $\mu_j>pL^{p-1}$, the argument in the proof of Theorem~\ref{main} shows that $v_j$ is unbounded, since that argument only uses $\mu_j>pL^{p-1}$. So we have that an index $j\geq2$ contributes to
the bounded kernel if and only if
\[
j(j+N-2)\leq pL^{p-1}.
\]
Solving this inequality for $j$ gives
$2\leq j\leq J(N,p)$, with $J(N,p)$ as in
\ref{eq:J-Np-JL}. The spherical-harmonic decomposition therefore gives \ref{kerfin}.\\
Arguing as in the proof of Theorem \ref{main} we deduce that all the newly constructed solutions \(v_j\) do not change sign for all $r>0$ and
    \begin{equation*}
v_j(r)\sim r^{\beta^+_j}\text{ for }r\to+\infty.
\end{equation*}
Indeed, if $v_j$ behaved like $r^{\beta^-_j}$, then $r^{N-3}v_0v_j$ would be integrable
and the boundary terms in \ref{a3}, by Lemma \ref{lem:asymptotics} (ii), would vanish as $R\to+\infty$, leading to the same
contradiction as in the proof of Theorem \ref{main}. 
For each $j\geq2$ 
the space of spherical harmonics of degree
$j$ has dimension
\[
\dim\mathcal H_j(\mathbb S^{N-1})
=\binom{N+j-1}{j}-\binom{N+j-3}{j-2}.
\]
The spherical-harmonic decomposition therefore gives
\ref{eq:kernel-dim-JL}.

We end the proof by showing the last assertion of the theorem.

Fix any integer $m\geq2$ and choose an integer $N$ such that
\begin{equation}\label{eq:choose-N-JL}
2+2(1+\sqrt2)m<N
<2+2(1+\sqrt2)(m+1).
\end{equation}
Such an integer exists because this open interval has length
$2(1+\sqrt2)>1$; it also satisfies $N\geq12$.
By the equivalence established in the preceding Proposition,
\ref{eq:choose-N-JL} implies
\[
\mu_m<\frac{(N-2)^2}{4}<\mu_{m+1}.
\]
In particular $\hat p_m(N)>p_{\mathrm{JL}}(N)$.
Choose any
\[
p\in\bigl(p_{\mathrm{JL}}(N),\hat p_m(N)\bigr].
\]
Then $\mu_m\leq pL^{p-1}$, while
\[
\mu_{m+1}>\frac{(N-2)^2}{4}>pL^{p-1}.
\]
Since $\mu_j$ increases strictly with $j$, precisely the degrees
$j=2,\ldots,m$ satisfy $\mu_j\leq pL^{p-1}$.
Hence $\lfloor J(N,p) \rfloor = m$.
\end{proof}

\begin{remark}\label{confronto}
    We observe that the additional bounded solutions constructed above do not belong to the weighted space considered by Dávila, del Pino, Musso and Wei in \cite{DavilaDelPinoMussoWei2007}. More precisely, recall that the relevant part of their weighted norm is
\[
\sup_{|x|\geq 1}|x|^\frac{2}{p-1}|\psi(x)|.
\]
By \ref{rp} for every \(j\geq2\), the corresponding nontrivial regular solution \(v_j\) of the linearized radial equation has infinite weighted norm, whenever \(p>p_{\mathrm{JL}}(N)\). Indeed, for all $j\geq 2$
\begin{equation*}
v_j(r)\sim r^{\beta^+_j}\text{ for }r\to+\infty.
\end{equation*}
and $\beta^+_j > -\frac{2}{p-1}$ for all $p>p_\mathrm{JL}$.
Thus, although these modes belong to \(L^\infty(\mathbb R^N)\), they do not belong to the weighted space used in \cite{DavilaDelPinoMussoWei2007}. This distinction explains why the existence of additional bounded kernel elements does not contradict the corresponding weighted linear theory.
\end{remark}

\section{Positivity of $v_0$}\label{s3}
As observed in the Introduction, the positivity of the function $v_0$ appears to be known among experts in the field. However, since we were unable to find a suitable reference, we have decided to include a proof here. We propose a blow-up argument, in order to show how the Lane–Emden equation is related to the Gelfand-type problem, namely
\begin{equation}\label{prob2}
    \begin{cases}
    -\Delta u_\la = \lambda (1+u_\la)^p & \text{ in } \Omega\subset\mathbb R^N\\
    u_\la>0 & \text{ in } \Omega \\
    u_\la = 0 & \text{ on } \partial \Omega
    \end{cases}
\end{equation}

It is well known that there exists a finite positive extremal parameter $\la^*$ such that \ref{prob2} has a minimal classical solution $u_\la\in C^0(\bar\Om)\cap C^2(\Om)$ if $0<\la<\la^*$, while no solution exists, even in the weak sense, for $\la>\la^*$. The set $\{u_\la: 0<\la<\la^*\}$ forms a branch of classical solutions increasing in $\la$. Its increasing pointwise limit $u^*(x):=\lim _{\la\to\la^*}u_\la(x)$ is a weak solution of \ref{prob2} for $\la=\la^*$, which is called the extremal solution of \ref{prob2}, see \cite{KK, KC, BCMR}.\\
The regularity and properties of extremal solutions depend strongly on the dimension $N$,
domain and the exponent $p$.
In the ball $B_1$, the following result holds (see \cite{JosephLundgren1973, BV})
\begin{thm}
    Suppose $N\geq11$, $p\geq p_\mathrm{JL}$ and $\Om=B_1$. Then $u^*$ is unbounded. In fact
\begin{equation}\label{las}
     \la^*=\frac{2}{p-1}\left(N-\frac{2p}{p-1}\right)
\end{equation}
    and
    $$u^*(x)=|x|^{-2/(p-1)}-1.$$
\end{thm}
\begin{remark}
    We note that
    $$\la^*=L^{p-1}.$$
\end{remark}
In the following, we use a standard scaling argument: since $u_\la$ is radial we will denote by $u_\la(r)$. By the previous Theorem $$\|u_\lambda\|_\infty = u_\lambda(0) \to u^*(0) = +\infty\quad\text{as}\quad \lambda \to \lambda^*.$$
Next we define
$$\tilde{u}_\lambda(r):=\frac{u_\lambda(\beta_\lambda r)}{ ||u_\la||_\infty}\quad\hbox{with }\beta_\lambda=\left(\frac{1}{||u_\la||_\infty}\right)^{(p-1)/2}.$$
A straightforward computation gives that $\tilde{u}_\lambda$
satisfies 
$$\begin{cases}
 -\tilde{u}''_\lambda(r)-\frac{N-1}{r}\tilde{u}'_\lambda(r)=\lambda\left(\frac{1} {||u_\la||_\infty}+\tilde{u}_\lambda(r)\right)^p &\text{in  } (0,||u_\la||_\infty^{(p-1)/2})\\
\tilde{u}_\lambda(0)=1
\end{cases}$$
Since $0\leq \tilde{u}_\lambda\leq 1$, by standard regularity estimates (see for example~\cite{GT}) there exists a radial function $u_\infty\in C^2(\R^N)$, $u_\infty\not\equiv0$ such that
\begin{equation}\label{u1}
\tilde{u}_\lambda\to u_\infty \hbox{ in }C^2_{loc}(\R^N)\hbox{ as }\lambda\to\lambda^*
\end{equation}
and $u_\infty$ solves
\begin{equation}\label{u2}
\begin{cases}
-\Delta u_\infty=\lambda^* u_\infty^p &\text{in  } \mathbb{R}^N\\
u_\infty(0)=1
\end{cases}
\end{equation}
\begin{lemma} \label{lemma}
    If $p\geq p_\mathrm{JL}$, we have that 
    $$z(r)=\frac{2}{p-1}u_\infty+ru_\infty'>0$$
    for all $r\geq0$.
\end{lemma}
The function $(\lambda^*)^{1/(p-1)}u_\infty$ is the radial solution of \ref{eq:LE} with value $(\lambda^*)^{1/(p-1)}$ at the origin. Since the positivity of $\frac{2}{p-1}u + ru'$ is invariant under the scaling $u \mapsto \sigma^{2/(p-1)}u(\sigma x)$, Lemma \ref{lemma} yields the positivity of $v_0$.
We have thus proved the following.
\begin{cor}\label{pos}
    If $p\geq p_\mathrm{JL}$, we have that 
    $$v_0(r)=\frac{2}{p-1}u+ru'>0$$
    for all $r\geq0$.
\end{cor}
We now prove Lemma \ref{lemma} 
\begin{proof}
   Suppose for the sake of contradiction that $z$ vanishes at some point $R > 0$ and $R$ is the first point where $z(r) = 0$. In the ball $B_R$, $z$ would satisfy:
   $$\begin{cases} 
   -\Delta z = \la^*p u_\infty^{p-1} z & \text{in } B_R \\ 
   z>0 & \text{in } B_R\\
   z= 0 & \text{on } \partial B_R 
   \end{cases}$$
   This implies that the first eigenvalue $\lambda_1(-\Delta - p \lambda^* u_{\infty}^{p-1}, B_R)$ must be equal to zero, since $z$ is the first eigenfunction of the operator $-\Delta - p \lambda^* u_{\infty}^{p-1}$.
   On the other hand for minimal solutions of \ref{prob2} we have that
   $$\la_1(-\Delta-p\la(1+u_\la)^{p-1})>0$$
   By \cite{BV}
   $$\lim_{\la\to\la^*}\la_1(-\Delta-\la p(1+u_\la)^{p-1})>0$$
   Scaling and passing to the limit, we obtain
   $$\la_1(-\Delta-p\la^*u^{p-1}_\infty, B_{2R})\geq\la_1(-\Delta-p\la^*u^{p-1}_\infty, \mathbb R^N)\geq 0$$
    By the strict monotonicity of the first Dirichlet eigenvalue with respect to the domain,
$$\lambda_1(-\Delta - p\lambda^* u_\infty^{p-1}, B_R) > \lambda_1(-\Delta - p\lambda^* u_\infty^{p-1}, B_{2R}) \ge 0.$$
This contradicts the previous finding that the eigenvalue was zero; hence, $z$ must remain strictly positive.
\end{proof}

\bibliographystyle{abbrv}
\bibliography{references_lane_emden}

@article{Aubin1976,
  author  = {Aubin, Thierry},
  title   = {Probl\`emes isop\'erim\'etriques et espaces de {Sobolev}},
  journal = {Journal of Differential Geometry},
  volume  = {11},
  number  = {4},
  pages   = {573--598},
  year    = {1976},
  doi     = {10.4310/JDG/1214433725}
}

@article{Talenti1976,
  author  = {Talenti, Giorgio},
  title   = {Best constant in {Sobolev} inequality},
  journal = {Annali di Matematica Pura ed Applicata},
  volume  = {110},
  pages   = {353--372},
  year    = {1976},
  doi     = {10.1007/BF02418013}
}

@article{CaffarelliGidasSpruck1989,
  author  = {Caffarelli, Luis A. and Gidas, Basilis and Spruck, Joel},
  title   = {Asymptotic symmetry and local behavior of semilinear elliptic equations with critical {Sobolev} growth},
  journal = {Communications on Pure and Applied Mathematics},
  volume  = {42},
  number  = {3},
  pages   = {271--297},
  year    = {1989},
  doi     = {10.1002/cpa.3160420304}
}

@article{JosephLundgren1973,
  author  = {Joseph, D. D. and Lundgren, T. S.},
  title   = {Quasilinear {Dirichlet} problems driven by positive sources},
  journal = {Archive for Rational Mechanics and Analysis},
  volume  = {49},
  pages   = {241--269},
  year    = {1973}
}

@article{GuiNiWang1992,
  author  = {Gui, Changfeng and Ni, Wei-Ming and Wang, Xuefeng},
  title   = {On the stability and instability of positive steady states of a semilinear heat equation in {$\mathbb{R}^n$}},
  journal = {Communications on Pure and Applied Mathematics},
  volume  = {45},
  number  = {9},
  pages   = {1153--1181},
  year    = {1992},
  doi     = {10.1002/cpa.3160450906}
}

@article{DavilaDelPinoMussoWei2007,
  author  = {D{\'a}vila, Juan and del Pino, Manuel and Musso, Monica and Wei, Juncheng},
  title   = {Standing waves for supercritical nonlinear {Schr\"odinger} equations},
  journal = {Journal of Differential Equations},
  volume  = {236},
  number  = {1},
  pages   = {164--198},
  year    = {2007},
  doi     = {10.1016/j.jde.2007.01.016}
}

@article{Farina2007,
  author  = {Farina, Alberto},
  title   = {On the classification of solutions of the {Lane--Emden} equation on unbounded domains of {$\mathbb{R}^N$}},
  journal = {Journal de Math\'ematiques Pures et Appliqu\'ees},
  volume  = {87},
  number  = {5},
  pages   = {537--561},
  year    = {2007},
  doi     = {10.1016/j.matpur.2007.03.001}
}

@article{DupaigneFarina2010,
  author  = {Dupaigne, Louis and Farina, Alberto},
  title   = {Stable solutions of {$-\Delta u=f(u)$} in {$\mathbb{R}^N$}},
  journal = {Journal of the European Mathematical Society},
  volume  = {12},
  number  = {4},
  pages   = {855--882},
  year    = {2010},
  doi     = {10.4171/JEMS/217}
}

@book{Dupaigne2011,
  author    = {Dupaigne, Louis},
  title     = {Stable Solutions of Elliptic Partial Differential Equations},
  series    = {Chapman \& Hall/CRC Monographs and Surveys in Pure and Applied Mathematics},
  volume    = {143},
  publisher = {Chapman \& Hall/CRC},
  year      = {2011}
}

@article{FrankKimPistoia2021,
  author  = {Frank, Rupert L. and Kim, Seunghyeok and Pistoia, Angela},
  title   = {Non-degeneracy for the critical {Lane--Emden} system},
  journal = {Proceedings of the American Mathematical Society},
  volume  = {149},
  number  = {1},
  pages   = {265--278},
  year    = {2021},
  doi     = {10.1090/proc/15217}
}

@article {BV,
    AUTHOR = {Brezis, Haim and V\'azquez, Juan Luis},
     TITLE = {Blow-up solutions of some nonlinear elliptic problems},
   JOURNAL = {Rev. Mat. Univ. Complut. Madrid},
  FJOURNAL = {Revista Matem\'atica de la Universidad Complutense de Madrid},
    VOLUME = {10},
      YEAR = {1997},
    NUMBER = {2},
     PAGES = {443--469},
      ISSN = {0214-3577},
   MRCLASS = {35J65 (35B05 35P30 35Q55)},
  MRNUMBER = {1605678},
MRREVIEWER = {C.\ A.\ Swanson},
}

@book {GT,
    AUTHOR = {Gilbarg, David and Trudinger, Neil S.},
     TITLE = {Elliptic partial differential equations of second order},
    SERIES = {Classics in Mathematics},
      NOTE = {Reprint of the 1998 edition},
 PUBLISHER = {Springer-Verlag, Berlin},
      YEAR = {2001},
     PAGES = {xiv+517},
      ISBN = {3-540-41160-7},
   MRCLASS = {35-02 (35Jxx)},
  MRNUMBER = {1814364},
}

@article {BCMR,
    AUTHOR = {Brezis, Ha\"im and Cazenave, Thierry and Martel, Yvan and
              Ramiandrisoa, Arthur},
     TITLE = {Blow up for {$u_t-\Delta u=g(u)$} revisited},
   JOURNAL = {Adv. Differential Equations},
  FJOURNAL = {Advances in Differential Equations},
    VOLUME = {1},
      YEAR = {1996},
    NUMBER = {1},
     PAGES = {73--90},
      ISSN = {1079-9389},
   MRCLASS = {35K60 (35J65)},
  MRNUMBER = {1357955},
MRREVIEWER = {Song\ Mu\ Zheng},
}

@article {KC,
    AUTHOR = {Keller, Herbert B. and Cohen, Donald S.},
     TITLE = {Some positone problems suggested by nonlinear heat generation},
   JOURNAL = {J. Math. Mech.},
  FJOURNAL = {J. Math. Mech.},
    VOLUME = {16},
      YEAR = {1967},
     PAGES = {1361--1376},
   MRCLASS = {35.10},
  MRNUMBER = {213694},
MRREVIEWER = {H.\ H.\ Schaefer},
}

@article {KK,
    AUTHOR = {Keener, J. P. and Keller, H. B.},
     TITLE = {Positive solutions of convex nonlinear eigenvalue problems},
   JOURNAL = {J. Differential Equations},
  FJOURNAL = {Journal of Differential Equations},
    VOLUME = {16},
      YEAR = {1974},
     PAGES = {103--125},
      ISSN = {0022-0396,1090-2732},
   MRCLASS = {35G20},
  MRNUMBER = {346305},
MRREVIEWER = {P.\ P.\ Zabre\u iko},
       DOI = {10.1016/0022-0396(74)90029-1},
       URL = {https://doi.org/10.1016/0022-0396(74)90029-1},
}

@article{Collot2017,
  author    = {Charles Collot},
  title     = {Nonradial type {II} blow up for the energy-supercritical semilinear heat equation},
  journal   = {Analysis \& {PDE}},
  volume    = {10},
  number    = {1},
  pages     = {127--252},
  year      = {2017},
  doi       = {10.2140/apde.2017.10.127},
}

@article{Rey1990,
  author    = {Olivier Rey},
  title     = {The role of the {Green}'s function in a nonlinear elliptic equation 
               involving the critical {Sobolev} exponent},
  journal   = {Journal of Functional Analysis},
  volume    = {89},
  number    = {1},
  pages     = {1--52},
  year      = {1990},
  doi       = {10.1016/0022-1236(90)90002-3},
}

@article {SJ,
    AUTHOR = {Serrin, James and Zou, Henghui},
     TITLE = {Classification of positive solutions of quasilinear elliptic
              equations},
   JOURNAL = {Topol. Methods Nonlinear Anal.},
  FJOURNAL = {Topological Methods in Nonlinear Analysis},
    VOLUME = {3},
      YEAR = {1994},
    NUMBER = {1},
     PAGES = {1--25},
      ISSN = {1230-3429},
   MRCLASS = {35J60 (34A34 35B05)},
  MRNUMBER = {1272885},
MRREVIEWER = {Allan\ L.\ Edelson},
       DOI = {10.12775/TMNA.1994.001},
       URL = {https://doi.org/10.12775/TMNA.1994.001},
}

@book{BoyceDiPrima2005,
  author    = {Boyce, William E. and DiPrima, Richard C.},
  title     = {Elementary Differential Equations and Boundary Value Problems},
  edition   = {8th},
  publisher = {John Wiley \& Sons},
  address   = {Hoboken, NJ},
  year      = {2005}
}

\end{document}